\documentclass[twocolumn,showkeys,preprintnumbers,amsmath,amssymb,superscriptaddress]{revtex4-1}

\usepackage{graphicx}
\usepackage{dcolumn}
\usepackage{bm}
\usepackage{soul, xcolor}
\usepackage{color, wasysym}
\usepackage{textcomp}
\usepackage{amsmath, amssymb, amsthm, latexsym, epsfig}
\usepackage{mathrsfs}

\DeclareMathOperator{\D}{d\!}

\def\ulamek#1#2{\mbox{\normalfont$\frac{#1}{#2}$}}

\begin{document}

\newtheorem{theorem}{Theorem}
\newtheorem{definition}{Definition}
\newtheorem{lemma}{Lemma}
\newtheorem{proposition}{Proposition}
\newtheorem{remark}{Remark}
\newtheorem{con}{Conjecture}
\newtheorem{example}{Example}


\title{Hausdorff moment problem for combinatorial numbers of Brown and Tutte: exact solution}

\author{K. A. Penson}
\email{karol.penson@sorbonne-universite.fr}
\affiliation{Laboratoire de Physique Th\'eorique de la Mati\`{e}re Condens\'{e}e (LPTMC), Sorbonne Universit\'{e}, Campus Pierre et Marie Curie (Paris 06), CNRS UMR 7600,\\ Tour 13 - 5i\`{e}me \'et., B.C. 121, 4 pl. Jussieu, F 75252 Paris Cedex 05 France}

\author{K. G\'{o}rska} 
\email{katarzyna.gorska@ifj.edu.pl}
\affiliation{Institute of Nuclear Physics, Polish Academy of Sciences, \\ ul. Radzikowskiego 152, PL-31342 Krak\'{o}w, Poland}
\affiliation{Laboratoire de Physique Th\'eorique de la Mati\`{e}re Condens\'{e}e (LPTMC), Sorbonne Universit\'{e}, Campus Pierre et Marie Curie (Paris 06), CNRS UMR 7600,\\ Tour 13 - 5i\`{e}me \'et., B.C. 121, 4 pl. Jussieu, F 75252 Paris Cedex 05 France}

\author{A. Horzela} 
\email{andrzej.horzela@ifj.edu.pl}
\affiliation{Institute of Nuclear Physics, Polish Academy of Sciences, \\ ul. Radzikowskiego 152, PL-31342 Krak\'{o}w, Poland}

\author{G. H. E. Duchamp}
\email{gheduchamp@gmail.com}
\affiliation{Laboratoire d'Informatique de Paris-Nord (LIPN), Sorbonne Universit\'{e}, Universit\'{e} Paris - Nord (Paris 13), CNRS UMR 7030,\\ Villetaneuse F 93430 France} 

\begin{abstract}
\medskip
We investigate the combinatorial sequences $A(M, n)$ introduced by W. G. Brown  (1964) and W. T. Tutte (1980) appearing in enumeration of convex polyhedra. Their formula is
\begin{equation*}
A(M, n) = \frac{2 (2M+3)!}{(M+2)! M!}\,\frac{(4n+2M+1)!}{n! (3n + 2M + 3)!}
\end{equation*}
with $n, M =0, 1, 2, \ldots$, and we conceive it as Hausdorff moments, where $M$ is a parameter and $n$ enumerates the moments. We solve exactly the corresponding Hausdorff moment problem: $A(M, n) = \int_{0}^{R} x^{n} W_{M}(x) \D x$ on the natural support $(0, R)$, $R = 4^{4}/3^{3}$, using the method of inverse Mellin transform. We provide explicitly the weight functions $W_{M}(x)$ in terms of the Meijer G-functions $G_{4, 4}^{4, 0}$, or equivalently, the generalized hypergeometric functions ${_{3}F_{2}}$ (for $M=0, 1$) and ${_{4}F_{3}}$ (for $M \geq 2$). For $M = 0, 1$, we prove that $W_{M}(x)$ are non-negative and normalizable, thus they are probability distributions. For $M \geq 2$, $W_{M}(x)$ are signed functions vanishing on the extremities of the support. By rephrasing this problem entirely in terms of Meijer G representations we reveal an integral relation which directly furnishes $W_M(x)$ based on ordinary generating function of $A(M, n)$ as an input. All the results are studied analytically as well as graphically.  
\end{abstract}

\maketitle

\section{Introduction}

Combinatorial numbers, which are by necessity positive integers, turn out very often to be related to probability, as they can be identified as power moments of positive and normalizable functions, i.e. the probability distributions. In most known cases the support of these distribution are either set of positive integers, the positive half-axis, or finite segments of the positive axis in the form $(0, R)$. For instance, many combinatorial numbers characterizing set partitions \cite{TM1, TM2} turn out to be moments of positive functions. Certain sequences of numbers contain parameters which permit to relate them to probability distributions for limited values of parameters only, see for instance \cite{WMlotkowski14}. (As we shall see later the sequence Eq. \eqref{22/06-6} also belongs to this category.) In this work we concentrate on a sequence of combinatorial numbers appearing in counting of bisections of convex polyhedra \cite{WGBrown64, WTTutte80}, which reads:
\begin{equation}\label{22/06-5}
A(M, n) = \frac{2 (2M+3)!}{(M+2)! M!}\,\frac{(4n+2M+1)!}{n! (3n + 2M + 3)!}, 
\end{equation}
where $M, n=0, 1, \ldots$. $A(M, n)$ are integers for all $M$ and $n$. We enumerate below the initial values $n=0, \ldots, 5$ of $A(M, n)$ for $0 \leq M \leq 4$:
\begin{align}\label{30/06-10}
A(0, n) &= 1, 1, 3, 13, 68, 399, 2530, \ldots \\ \label{30/06-11}
A(1, n) &= 2, 5, 20, 100, 570, 3542, \ldots \\ \label{30/06-12}
A(2, n) &= 5, 21, 105, 595, 3675, 24150, \ldots \\ \label{30/06-13}
A(3, n) &= 14, 84, 504, 3192, 21252, 147420, \ldots \\ \label{30/06-14}
A(4, n) &= 42, 330, 2310, 16170, 115500, 844074, \ldots\!.
\end{align}
The sequences $A(M, n)$ for $0 \leq M \leq 3$ are documented and discussed in N. J. A. Sloane's Online Encyclopedia of Integer Seqences (OEIS) \cite{NJAS}: \\ $A(0, n) = $A000260(n), $A(1, n) =$ A197271(n), $A(2, n) =$ A341853(n), and $A(3, n) =$ A341854(n).
However notice that
\begin{equation}\label{28/06-5}
A(M, 0) = {\rm Cat}(M+1), 
\end{equation}
where ${\rm Cat}(n) = \binom{2n}{n} \frac{1}{n+1}$ are Catalan numbers. We set out to solve the following Hausdorff power  moment problem: find $W_{M}(x)$ satisfying the infinite set of equations:
\begin{equation}\label{22/06-6}
A(M, n) = \int_{0}^{R} x^{n}\, W_{M}(x) \D x, \qquad n=0, 1, \ldots,
\end{equation}
where $R$ is given by the known formula $R = \lim_{n\to\infty} [A(M, n)]^{1/n} = 4^{4}/3^{3}$, i.e. is independent on $M$. We shall employ the method of inverse Mellin transform, which implies for $n=s - 1$ that
\begin{equation}\label{22/06-7}
W_{M}(x) = \mathcal{M}^{-1}[A(M, s-1); x].
\end{equation}

In Sec. \ref{sec2} we shall enumerate and present the conventional tools applied in calculating the inverse Mellin transforms, including the Meijer G-functions, the generalized hypergeometric functions and some of their properties. With the above tools at hand, in Sec. \ref{sec3} we shall perform in detail the Mellin inversion and obtain the explicit closed-form solutions for $W_{M}(x)$ in terms of generalized hypergeometric functions ${_{3}F_{2}}$ and ${_{4}F_{3}}$. We prove the positivity of $W_{M}(x)$ for $M=0, 1$ only, whereas for $M \geq 2$ we demonstrate that $W_{M}(x)$ are signed functions. $W_{M}(x)$ are discussed graphically for $0 \leq M \leq 4$. In Sec. \ref{sec4} we derive the closed-form expression for the ordinary generating function (ogf) of $A(M, n)$, i.e. for $G(M, z) = \sum_{n=0}^{\infty} A(M, n) z^{n}$, and we establish a relationship between $W_{M}(x)$ and $G(M, z)$, by rephrasing both of them in the common language of Meijer G-functions. The above mentioned relation allows to construct explicitly the function $W_{M}(x)$ using {\em solely} the hypergeometric representation of $G(M, z)$. This procedure is strongly evocative of the inversion of one-sided, finite Hilbert transform. Sec. \ref{sec5} contains the discussion and final remarks.
 
\section{Definitions and preliminaries}\label{sec2}

The main tool that we employ in the treatment of various moment problems is the Mellin transform $\mathcal{M}$ and it inverse $\mathcal{M}^{-1}$. In the following we group several definitions and informations about the Mellin transform of a function $f(x)$ defined for $x\geq 0$. The Mellin transform is defined for complex $s$ as \cite{INSneddon72}
\begin{equation}\label{27.10.13-1}
\mathcal{M}[f(x); s] = f^{\star}(s) = \int_{0}^{\infty} x^{s-1} f(x) dx,
\end{equation}
along with its inverse
\begin{equation}\label{27.10.14-2}
\mathcal{M}^{-1}[f^{\star}(s); x] = f(x) = \frac{1}{2\pi i} \int_{c-i\infty}^{c+i\infty} x^{-s} f^{\star}(s) ds.
\end{equation}
For the role of constant $c$ consult \cite{INSneddon72}.

If $\mathcal{M}[f(x); s] = f^{\star}(s)$ and $\mathcal{M}[g(x); s] = g^{\star}(s)$ then
\begin{align}\label{16/06-1}
\begin{split}
\mathcal{M}^{-1}[f^{\star}(s) g^{\star}(s); x]  & = \int_{0}^{\infty} f\left(\frac{x}{t}\right) g(t) \frac{1}{t} \D t \\
& =  \int_{0}^{\infty} g\left(\frac{x}{t}\right) f(t) \frac{1}{t} \D t.
\end{split}
\end{align}
The last two integrals are called Mellin (i.e. multiplicative) convolutions of $f(x)$ with $g(x)$. For fixed $a>0$, $h\neq 0$, the Mellin transform satisfies the following scaling property:
\begin{equation}\label{21/06-4}
\mathcal{M}[x^{b}f(ax^{h}); s] = \frac{1}{|h|} a^{-\frac{s+b}{h}} f^{\star}(\ulamek{s+b}{h}).
\end{equation}

Among known Mellin transforms a very special role is played by those entirely expressible through products and ratios of Euler's gamma functions. The Meijer G-function is defined as an inverse Mellin transform \cite{APPrudnikov-v3}:
\begin{align}\label{16/06-2}
& G^{m, n}_{p, q}\left(x\Big\vert{\alpha_{1} \ldots \alpha_{p} \atop \beta_{1} \ldots \beta_{q}}\right) \nonumber\\ &= \mathcal{M}^{-1}\left[\frac{\prod_{j=1}^{m}\Gamma(\beta_{j}+s)\, \prod_{j=1}^{n}\Gamma(1-\alpha_{j}-s)}{\prod_{j=m+1}^{q}\Gamma(1-\beta_{j}-s)\, \prod_{j=n+1}^{p}\Gamma(\alpha_{j}+s)}; x\right] \\ 
&\qquad = {\rm MeijerG}([[\alpha_{1}, \ldots, \alpha_{n}], [\alpha_{n+1}, \ldots, \alpha_{p}]], \nonumber \\ & \qquad\qquad\qquad\quad [[\beta_{1}, \ldots, \beta_{m}], [\beta_{m+1}, \ldots, \beta_{q}]], x). \label{15/06-10}
\end{align}
The notation for $G^{m, n}_{p, q}$ in Eq. \eqref{15/06-10} is motivated by Maple and Mathematica notation \cite{F1}. We will consequently use both notations throughout this paper. In Eq. \eqref{16/06-2} empty products are taken to be equal to 1. In Eqs. \eqref{16/06-2} and \eqref{15/06-10} the parameters are subject of conditions:
\begin{align}\label{21/06-1}
\begin{split}
&z\neq 0, \quad 0 \leq m \leq q, \quad 0 \leq n \leq p, \\
&\alpha_{j} \in \mathbb{C}, \quad j=1,\ldots, p; \quad \beta_{j}\in\mathbb{C}, \quad j=1,\ldots, q. 
\end{split}
\end{align}
See Refs. \cite{APPrudnikov-v3, KGorska13} for a full description of integration contours in Eq. \eqref{16/06-2}, general properties and special cases of the Meijer G-functions. The convergence of the Mellin inversion in Eqs. \eqref{16/06-2} and \eqref{15/06-10} is conditioned upon specific requirements involving both chains of parameters $(a_{p})$ and $(b_{q})$. The aforementioned conditions will be quoted and checked in Sec. \ref{sec3} on the example of the Mellin inversion derived from Eq. \eqref{22/06-5}.

The generalized hypergeometric function $_{p}F_{q}$ is defined as:
\begin{multline}\label{15/06-11}
{_{p}F_{q}}\left({a_{1}, \ldots, a_{p} \atop b_{1}, \ldots, b_{q}}; z\right) =  \sum_{k=0}^{\infty} \frac{(a_{1})_{k} \cdots (a_{p})_{k}}{(b_{1})_{k} \cdots (b_{q})_{k}}\, \frac{z^{k}}{k!} \\
=  {_{p}F_{q}}([a_{1}, \ldots, a_{p}], [b_{1}, \ldots, b_{q}]; z),
\end{multline}
where $(a)_{k} = \Gamma(a + k)/\Gamma(a)$ is called the Pochhammer symbol, and neither of $b_{j}$, $j=1, \ldots, q$, is a negative integer, see \cite{JThomae1870}. 

We shall also use the following relation linking one ${_{p}F_{q}}$ with one $G_{p, q}^{m, n}$, for $p \leq q+1$:
\begin{multline}\label{15/06-12}
{_{p}F_{q}}\left({a_{1},\ldots, a_{p} \atop b_{1}, \ldots, b_{q}}; z\right) = \left(\prod_{k=1}^{q}\Gamma(b_{k})/\prod_{k=1}^{p}\Gamma(a_{k})\right)\\ \times G^{1, p}_{p, q+1}\left(-z\Big\vert {1-a_{1}, \ldots, 1 - a_{p} \atop 0, 1-b_{1}, \ldots, 1-b_{q}} \right)\!,
\end{multline}
see Eq. (16.18.1) of \cite{NIST}, where a particular attention should be paid to the position of $0$ in the lower list of parameters. For the proof of Eq. \eqref{15/06-12} see Eq. (12.3.18) on p. 317 of \cite{RBeals16}. Additional identities satisfied by $G^{m, n}_{p, q}$ and used in this work are
\begin{align}\label{15/06-13}
G^{m, n}_{p, q}\left(\frac{1}{z} \Big\vert {a_{1}, \ldots, a_{p} \atop b_{1}, \ldots, b_{q}} \right) & = G^{n, m}_{q, p}\left(z \Big\vert {1-b_{1}, \ldots, 1-b_{q} \atop 1-a_{1}, \ldots, 1-a_{p}} \right) \\
z^{\mu}\, G^{m, n}_{p, q}\left(z \Big\vert {a_{1}, \ldots, a_{p} \atop b_{1}, \ldots, b_{q}} \right) & =  G^{m, n}_{p, q}\left(z \Big\vert {a_{1} + \mu, \ldots, a_{p} + \mu \atop b_{1} + \mu, \ldots, b_{q} + \mu} \right), \label{15/06-13b}
\end{align}
see Eqs. (16.19.1) and (16.19.2) of \cite{NIST}, correspondingly. 

For certain conditions satisfied by the parameter lists, the functions $G^{m, n}_{p, q}$ can be represented as a finite sum of hypergeometric function, see Eq. (16.17.2) of \cite{NIST} and/or Eq. (8.2.2.3) of \cite{APPrudnikov-v3}, which are sometimes referred to as Slater relations.

We quote for reference the Gauss-Legendre multiplication formula for gamma function encountered in this work:
\begin{align}\label{27.10.13-9}
\begin{split}
&\Gamma(nz) = (2\pi)^{\frac{1-n}{2}} n^{nz-\frac{1}{2}} \prod_{j=0}^{n-1}\Gamma\left(z+\frac{j}{n}\right), \\
&z\neq 0, -1, -2, \ldots, \quad n=1, 2, \ldots. 
\end{split}
\end{align}
We also introduce a short notation for a special list of $k$ elements:
\begin{equation}\label{27.10.13-10}
\Delta(k, a) = \frac{a}{k}, \frac{a+1}{k}, \ldots, \frac{a + k-1}{k}, \qquad k\neq 0.
\end{equation}

\section{Solving the moment problem}\label{sec3}

In this section we shall derive the exact and explicit forms of the solutions $W_{M}(x)$ of the Hausdorff moment problem \eqref{22/06-6} where $A(M, n)$ is given by Eq. \eqref{22/06-5}. Denote 
\begin{equation}\label{28/06-6}
P(M)= \frac{2(2M + 3)!}{(M+2)! M!}, 
\end{equation}
set in Eq. \eqref{22/06-6} $n = s - 1$, and use twice the Gauss-Legendre formula Eq. \eqref{27.10.13-9} in transforming Eq. \eqref{22/06-6}  to obtain $A(M, s-1) \equiv \tilde{A}(M, s)$:
\begin{multline}\label{23/06-1}
\tilde{A}(M, s) = r_{W}(M)R^{\,s}\, \\ \times \frac{\Gamma(s-\frac{1}{2} + \frac{M}{2}) \Gamma(s-\frac{1}{4} + \frac{M}{2}) \Gamma(s + \frac{M}{2}) \Gamma(s + \frac{1}{4} + \frac{M}{2})}{\Gamma(s + \frac{1}{3} + \frac{2M}{3}) \Gamma(s + \frac{2}{3} + \frac{2M}{3}) \Gamma(s + 1 + \frac{2M}{3}) \Gamma(s)},
\end{multline}
where
\begin{equation}
r_{W}(M) = \frac{3^{\frac{1}{2} - 2M} 2^{4M + \frac{1}{2}}}{192 \sqrt{\pi}} P(M).
\end{equation}
We apply now the scaling property Eq. \eqref{21/06-4} along with the definitions of the Meijer G-function Eqs. \eqref{16/06-2} and \eqref{15/06-10} in order to write the final form of $W_{M}(x) = \mathcal{M}^{-1}[\tilde{A}(M, s); x]$:
\begin{align}\label{23/06-5}
&W_{M}(x) =  r_{W}(M)\, G_{4, 4}^{\,4, 0}\left(\frac{x}{R}\Big\vert {0, \Delta(3, 2M+1) \atop \Delta(4, 2M-2)}\right) \quad\\
&= r_{W}(M)\, {\rm MeijerG}\Big(\big[[\,\,\,], [0, \ulamek{2M}{3} + \ulamek{1}{3},  \ulamek{2M}{3} + \ulamek{2}{3}, \ulamek{2M}{3} + 1]\big], \nonumber \\ &\qquad\qquad \left. \big[[\ulamek{M}{2} - \ulamek{1}{2}, \ulamek{M}{2} - \ulamek{1}{4}, \ulamek{M}{2}, \ulamek{M}{2} + \ulamek{1}{4}], [\,\,\,]\big], \frac{x}{R}\right). \label{23/06-6}
\end{align}
The solutions in Eq. \eqref{23/06-5} are unique. According to the definitions of Eqs. \eqref{16/06-2} and \eqref{15/06-10} the parameter lists $(\alpha_{1}, \ldots, \alpha_{4})$ and $(\beta_{1}, \ldots, \beta_{4})$ for $p=q=4$ in Eq. \eqref{23/06-5} can be read off as $(\alpha_{p}) = (0, \Delta(3, 2M+1)\,)$ and $(\beta_{q}) = (\,\Delta(4, 2M-2)\,)$, using Eq. \eqref{27.10.13-10}. We can now extract the conditions for convergence of integral \eqref{16/06-2} as a function of $(\alpha_{p})$ and $(\beta_{q})$. They define the range of variable $s$ for which the convergence is assured with the formula (2.24.2.1) of \cite{APPrudnikov-v3}. Here for $m=4$, $n=0$, and $p = q = 4$, and the auxiliary parameter $c^{\star} \equiv m+n - (p+q)/2 = 0$. Thus, the range of real $s$ is determined from the inequality:
\begin{align}\label{28/06-1}
- \min_{1 \leq j \leq m}(\beta_{j}) \leq &s \leq 1 - \max_{1 \leq j \leq n}(\alpha_{j}), \quad \text{which reads} \\
\frac{1}{4} - \frac{M}{2} \leq &s \leq 1 -(-\infty), \qquad \text{or for} \quad s=n'+1 \nonumber\\
\frac{1}{4} - \frac{M}{2} \leq &n'+1 \leq \infty, \qquad \text{and finally} \nonumber \\
 -\frac{3}{4}-\frac{M}{2} \leq &n' \leq \infty, \label{28/06-1a}
\end{align}
where $n'$ enumerates the moments. We conclude that for $W_{M}(x)$ all the moments $\int_{0}^{R} x^{n'} W_{M}(x) \D x$, for $0 \leq n' < \infty$ are legitimate and converging.

Before embarking on detailed evaluation of Eq. \eqref{23/06-5} we claim that for $M=0, 1$ the weight function $W_{M}(x)$ will be a positive function on $x\in(0, R)$. This is based on the Mellin convolution property of Eq. \eqref{16/06-1} which shows that if two individual functions are positive, then for positive arguments, their Mellin convolution is also positive. The second element of this reasoning tell us that 
\begin{multline}\label{26/06-1}
\mathcal{M}^{-1}\left[\frac{\Gamma(s+a)}{\Gamma(s + b)}; x\right] = \frac{(1-x)^{1-a+b} x^{a}}{\Gamma(b-a)} > 0, \quad \text{for} \\ 0 < x < 1, \,\,\, b > a, 
\end{multline}
which is the direct consequence of Eq. (8.4.2.3) on p. 631 of \cite{APPrudnikov-v3}.  Eq. \eqref{26/06-1} is strongly reminiscent of the classical Euler Beta function. Moreover, the Beta distribution is the probability measure characterised by the density function $g_{\alpha, \beta}(x) = \{\Gamma(\alpha+\beta)/[\Gamma(\alpha)\Gamma(\beta)]\}\, x^{\alpha-1} (1-x)^{\beta - 1}$ \cite{NBalakrishnan03}. The r.h.s. of Eq. \eqref{26/06-1} for $0 < x < 1$ and $b > a$ is a positive function.

Suppose that we will be able to order the shifts in four gamma ratios in Eq. \eqref{23/06-1} in such a way that for every ratio $b > a$, as in Eq. \eqref{26/06-1}. Then the resulting weight function will be a {\em threefold} Mellin convolution of positive functions, and, through the above argument, will itself be positive. Let us first enumerate the gamma shifts for $M = 0$ in Eq. \eqref{23/06-5}, with $u =$ upper and $l =$ lower shifts. In the formulas \eqref{26/06-2}, \eqref{26/06-4}, and \eqref{26/06-6} below, the arrow $"\Longrightarrow"$ should be understood as: "can be reordered as".
\begin{equation}\label{26/06-2}
M=0: \quad \left\{ {u: 0, 1, \frac{2}{3}, \frac{1}{3} \atop l: 0, -\frac{1}{2}, \frac{1}{4}, -\frac{1}{4}} \right\} \Longrightarrow \left\{ {u: 0, \frac{1}{3},  \frac{2}{3}, 1 \atop l: -\frac{1}{4}, -\frac{1}{2}, 0, \frac{1}{4}} \right\}
\end{equation}
resulting in the gamma ratios:
\begin{equation}\label{26/06-3}
\frac{\Gamma(s - \frac{1}{4})}{\Gamma(s + 0)}\, \frac{\Gamma(s - \frac{1}{2})}{\Gamma(s + \frac{1}{3})}\, \frac{\Gamma(s + 0)}{\Gamma(s + \frac{2}{3})}\, \frac{\Gamma(s + \frac{1}{4})}{\Gamma(s + 1)}.
\end{equation}
Then the resulting $W_{0}(x)$ will be a positive function. We continue with the same argument for $M = 1$:
\begin{equation}\label{26/06-4}
M=1: \quad \left\{ {u: 0, \frac{5}{3}, \frac{4}{3}, 1 \atop l: \frac{1}{2}, 0, \frac{3}{4}, -\frac{1}{4}} \right\} \Longrightarrow \left\{ {u: 0, 1, \frac{4}{3},  \frac{5}{3} \atop l: -\frac{1}{4}, 0, \frac{1}{2}, \frac{3}{4}} \right\},
\end{equation}
resulting in the gamma ratios:
\begin{equation}\label{26/06-5}
\frac{\Gamma(s - \frac{1}{4})}{\Gamma(s + 0)}\, \frac{\Gamma(s + 0)}{\Gamma(s + 1)}\, \frac{\Gamma(s + \frac{1}{2})}{\Gamma(s + \frac{4}{3})}\, \frac{\Gamma(s + \frac{3}{4})}{\Gamma(s + \frac{5}{3})}.
\end{equation}
Then again, the resulting $W_{1}(x)$ will be a positive function. The situation changes for $M=2$, as then
\begin{equation}\label{26/06-6}
M=2: \quad \left\{ {u: 0, \frac{7}{3}, 2, \frac{5}{3} \atop l: 1, \frac{1}{2}, \frac{5}{4}, \frac{3}{4}} \right\} \Longrightarrow \left\{ {u: 0, \frac{5}{3}, 2, \frac{7}{3} \atop l: \frac{1}{2}, 1, \frac{3}{4}, \frac{5}{4}} \right\},
\end{equation}
and the resulting gamma ratio
\begin{equation}\label{26/06-7}
\frac{\Gamma(s + \frac{1}{2})}{\Gamma(s + 0)}\, \frac{\Gamma(s + 1)}{\Gamma(s + \frac{5}{3})}\, \frac{\Gamma(s + \frac{3}{4})}{\Gamma(s + 2)}\, \frac{\Gamma(s + \frac{5}{4})}{\Gamma(s + \frac{7}{3})}
\end{equation}
excludes the positivity of $W_{2}(x)$ as here $$\mathcal{M}^{-1}\left[\Gamma\Big(s + \frac{1}{2}\Big)/\Gamma(s); x\right] \sim - \sqrt{x}/(1 - x)^{3/2} < 0$$ for $0 < x < 1$. Similar arguments exclude the positivity for $M \geq 2$. The method of studying the positivity via multiple Mellin convolution was initiated in \cite{KAPenson11} and further applied in \cite{ABostan, WMlotkowski13, KGorska13, WMlotkowski13a, KGorska13a}, to various sequences of combinatorial numbers.

Since $A(M, 0) \neq 1$, see Eq. \eqref{28/06-5}, it is reasonable not to compare $W_{M}(x)$ for different $M$, but rather to consider $\tilde{W}_{M}(x) = W_{M}(x)/A(M, 0)$, "normalized" weight functions for different $M$. Note that zeroth moments of $\tilde{W}_{M}(x)$ are equal to $1$, but higher moments of $\tilde{W}_{M}(x)$ are not anymore integers but are rationals. In order to do so, we represent the Meijer G-functions of Eqs. \eqref{23/06-5} as a finite sum of three generalized hypergeometric functions ${_{3}F_{2}}$ (for $M=0, 1$), and ${_{4}F_{3}}$ (for $M \geq 2$), employing Eq. (8.2.2.3) of \cite{APPrudnikov-v3}. This last formula also permits to write down the general expression for $\tilde{W}_{M}(x) = W_{M}(x)/A(M, 0)$ for arbitrary integer $M$ with the help of generalized hypergeometric functions. However, due to its complexity we shall not reproduce it here. Instead we quote below the explicit forms for $\tilde{W}_{M}(x)$ for $0 \leq M \leq 3$, with $R = 4^{4}/3^{3}$:
\begin{align}\label{24/06-1}
\tilde{W}_{0}(x) &= \frac{2}{\pi \sqrt{x}}\, {_{3}F_{2}}\left({-\frac{1}{2}, -\frac{1}{6}, \frac{1}{6} \atop \frac{1}{4}, \frac{3}{4}}, \frac{x}{R} \right) \nonumber\\ &- \frac{\sqrt{2}}{\pi x^{1/4}}\, {_{3}F_{2}}\left({-\frac{1}{4}, \frac{1}{12}, \frac{5}{12} \atop \frac{1}{2}, \frac{5}{4}}, \frac{x}{R} \right) \nonumber\\
& + \frac{\sqrt{2}\, x^{1/4}}{32 \pi}\, {_{3}F_{2}}\left({\frac{1}{4}, \frac{7}{12}, \frac{11}{12} \atop \frac{3}{2}, \frac{7}{4}}, \frac{x}{R} \right) \\
&= \frac{2}{\pi \sqrt{x}}\, {_{3}F_{2}}\big(\big[\!-\ulamek{1}{2}, -\ulamek{1}{6}, \ulamek{1}{6}\big], \big[\ulamek{1}{4}, \ulamek{3}{4}\big], \ulamek{x}{R}\big) \nonumber\\ &- \frac{\sqrt{2}}{\pi x^{1/4}}\, {_{3}F_{2}}\big(\big[\!-\ulamek{1}{4}, \ulamek{1}{12}, \ulamek{5}{12}\big], \big[\ulamek{1}{2}, \ulamek{5}{4}\big], \ulamek{x}{R}\big) \nonumber\\
& + \frac{\sqrt{2}\, x^{1/4}}{32 \pi}\, {_{3}F_{2}}\big(\big[\ulamek{1}{4}, \ulamek{7}{12}, \ulamek{11}{12}\big], \big[\ulamek{3}{2}, \ulamek{7}{4}\big], \ulamek{x}{R}\big). \label{24/06-1a}
\end{align}
In the following three equations we skip the Maple notation.
\begin{align}\label{24/06-2}
\tilde{W}_{1}(x) &= \frac{2\sqrt{2}}{\pi} x^{1/4}\, {_{3}F_{2}}\left({-\frac{5}{12}, -\frac{1}{12}, \frac{1}{4} \atop \frac{1}{2}, \frac{3}{4}}, \frac{x}{R} \right) \nonumber\\ 
&- \frac{5\sqrt{x}}{2\pi}\, {_{3}F_{2}}\left({-\frac{1}{6}, \frac{1}{6}, \frac{1}{2} \atop \frac{3}{4}, \frac{5}{4}}, \frac{x}{R} \right) \nonumber\\
& +  \frac{5\sqrt{2}}{16\pi} x^{3/4}\, {_{3}F_{2}}\left({\frac{1}{12}, \frac{5}{12}, \frac{3}{4} \atop \frac{5}{4}, \frac{3}{2}}, \frac{x}{R} \right),  
\end{align}
\begin{align}\label{24/06-3}
\tilde{W}_{2}(x) &= -\frac{14\sqrt{x}}{5\pi}\, {_{4}F_{3}}\left({-\frac{5}{6}, -\frac{1}{2}, -\frac{1}{6}, \frac{3}{2} \atop \frac{1}{4}, \frac{1}{2}, \frac{3}{4}}, \frac{x}{R} \right) \nonumber \\
& + \frac{3\sqrt{2}}{\pi} x^{3/4}\, {_{4}F_{3}}\left({-\frac{7}{12}, -\frac{1}{4}, \frac{1}{12}, \frac{7}{4} \atop \frac{1}{2}, \frac{3}{4}, \frac{5}{4}}, \frac{x}{R} \right) \nonumber\\
& - \frac{35\sqrt{2}}{32\pi} x^{5/4}\, {_{4}F_{3}}\left({-\frac{1}{12}, \frac{1}{4}, \frac{7}{12}, \frac{9}{4} \atop \frac{5}{4}, \frac{3}{2}, \frac{7}{4}}, \frac{x}{R} \right),  
\end{align}
and
\begin{align}\label{24/06-4}
\tilde{W}_{3}(x) &= -\frac{4\sqrt{2}}{\pi} x^{5/4}\, {_{4}F_{3}}\left({-\frac{3}{4}, -\frac{5}{12}, -\frac{1}{12}, \frac{9}{4} \atop \frac{1}{2}, \frac{3}{4}, \frac{5}{4}}, \frac{x}{R} \right) \nonumber\\
& + \frac{9}{\pi} x^{3/2}\, {_{4}F_{3}}\left({-\frac{1}{2}, -\frac{1}{6}, \frac{1}{6}, \frac{5}{2} \atop \frac{3}{4}, \frac{5}{4}, \frac{3}{2}}, \frac{x}{R} \right) \nonumber\\
& - \frac{21\sqrt{2}}{8\pi} x^{7/4}\, {_{4}F_{3}}\left({-\frac{1}{4}, \frac{1}{12}, \frac{5}{12}, \frac{11}{4} \atop \frac{5}{4}, \frac{3}{2}, \frac{7}{4}}, \frac{x}{R} \right).  
\end{align}

\begin{figure}[!h]
\begin{center}
\includegraphics[scale=0.2]{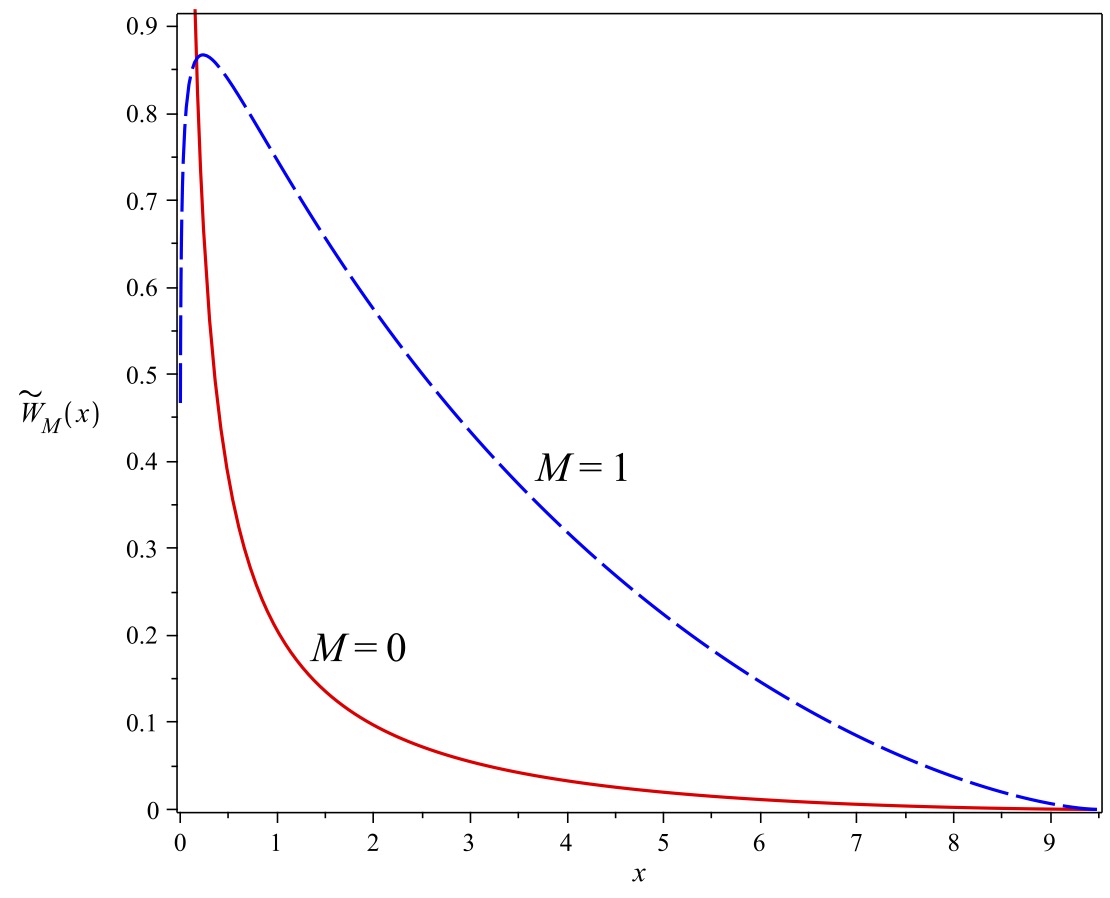}
\caption{\label{fig1}(Color online) Plot of $\tilde{W}_{M}(x)$ for $M=0$ (red continuous curve) and $M=1$ (blue dashed curve) for $x\in (0, R)$. Notice that $\tilde{W}_{0}(x)$ tends to infinity at $x=0$ whereas $\tilde{W}_{1}(x)$ approaches zero at $x=0$. $\tilde{W}_{0}(x)$ and $\tilde{W}_{1}(x)$ are normalized probability distributions.
} 
\includegraphics[scale=0.2]{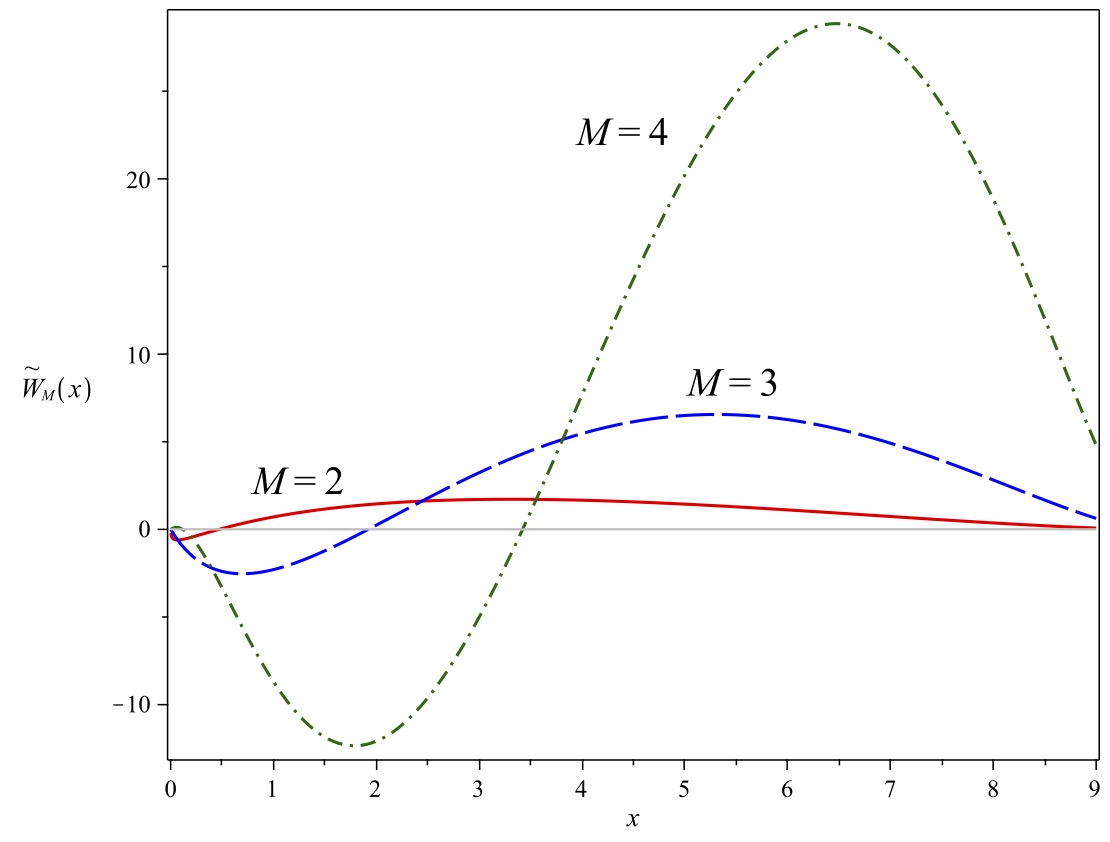}
\caption{\label{fig2}(Color online) Plot of $\tilde{W}_{M}(x)$ for $M=2$ (red continuous curve), $M=3$ (blue dashed curve), and $M=4$ (green dashed-dotted curve) for $x\in (0, R-0.48)$. Notice that $\tilde{W}_{M}(x)$ for $M \geq 2$ have a negative part and tend to zero at $x=0$.
} 
\end{center}
\end{figure}

We represent graphically $\tilde{W}_{M}(x)$ for $M = 0, 1$ on Fig. \ref{fig1}, and for $M = 2, 3, 4$ on Fig. \ref{fig2}.

\section{Linking generating and weight functions}\label{sec4}

We start with a general Hausdorff moment problem in the form of Eq. \eqref{22/06-6}. Solving Eq. \eqref{22/06-6} means to obtain $W(x)$ given the set $\rho(n)$, $n=0, 1, \ldots$. We define the ordinary generating function (ogf) of moments $\rho(n)$ as
\begin{equation}\label{30/06-1}
G(z) = \sum_{n=0}^{\infty} \rho(n) z^{n},
\end{equation}
with the radius of convergence equal to $1/R$, i.e. $z < 1/R$. We classically observe that
\begin{align}\label{30/06-2}
\begin{split}
G(z) & = \sum_{n=0}^{\infty} z^{n} \left[\int_{0}^{R} x^{n} W(x) \D x\right] \\ 
& = \int_{0}^{R} W(x) \left[\sum_{n=0}^{\infty} (x z)^{n}\right] \D x \\
& = \int_{0}^{R} \frac{W(x)}{1 - z x} \D x,
\end{split}
\end{align}
with $z x < 1$. Since in Eq. \eqref{30/06-2} $0 \leq x \leq R$, it implies $z < 1/R$. From Eq. \eqref{30/06-2} it follows that
\begin{equation}\label{30/06-3}
\frac{1}{z} G\left(\frac{1}{z}\right) = \int_{0}^{R} \frac{W(x)}{z - x} \D x, \qquad \text{with} \quad z > R.
\end{equation}
The above equations constitute the seed of the inversion procedure by Stieltjes to solve Eq. \eqref{30/06-3} using the complex analysis. For singularity analysis of $G(z)$ in complex plane see \cite{JGLiu16}. For a recent detailed application of this method, along with the exhaustive reference list, see \cite{ABostan20}. A very complete exposition of the Stieltjes method can be found in \cite{AHora07}. 

The above transformations are fairly standard, however in view of results of Sec. \ref{sec3}, a certain pattern does appear that permits one to deduce $W(x)$ directly from $G(z)$, via Eq. \eqref{30/06-3}. In order to make explicit this pattern several manipulations with $G(z) \equiv G(M, z)$ are needed.

We use the definition of Pochhammer symbols to write down the ogf $G(M, z)$ of the moments $A(M, n)$, and it reads
\begin{align}\label{30/06-4}
&G(M, z) = \frac{2 (2M+1)!}{(M+2)! M!}\, {_{4}F_{3}}\left({\Delta(4, 2M + 2) \atop \Delta(3, 2M + 4)}; R z \right)\\
 &= \frac{2 (2M+1)!}{(M+2)! M!}\, {_{4}F_{3}}\big(\big[1 + \ulamek{M}{2}, \ulamek{3}{4} + \ulamek{M}{2}, \ulamek{1}{2} + \ulamek{M}{2}, \ulamek{5}{4} + \ulamek{M}{2}\big], \nonumber\\ & \qquad\qquad\quad \big[2 + \ulamek{2M}{3}, \ulamek{5}{3} + \ulamek{2M}{3}, \ulamek{4}{3} +\ulamek{2M}{3}\big]; R z \big), \quad z< R. \label{30/06-4a}
\end{align}
We come back to Eq. \eqref{15/06-12} in order to frame Eq. \eqref{30/06-4} in the Meijer G-notation. Carrying out the products of gamma functions in Eq. \eqref{15/06-12} this furnishes:
\begin{align}\label{1/07-1}
& G(M, z) = r_{G}(M) G^{1, 4}_{4, 4}\left(-Rz \Big\vert {\Delta(4, -2M - 1) \atop 0, \Delta(3, -2M-3)}\right) \\
&= r_{G}(M) {\rm MeijerG}\big(\big[\big[-\ulamek{M}{2}, \ulamek{1}{2}-\ulamek{M}{2},\ulamek{1}{4}-\ulamek{M}{2}, -\ulamek{1}{4}-\ulamek{M}{2}\big], \big[\,\,\,\big]\big], \nonumber\\ &\qquad\qquad \big[\big[0\big], \big[-\ulamek{2}{3}-\ulamek{2M}{3}, -\ulamek{1}{3}-\ulamek{2M}{3}, -1-\ulamek{2M}{3} \big]\big], -Rz \big)
\end{align}
with
\begin{equation}\label{1/07-2}
r_{G}(M) = \frac{4}{81\sqrt{\pi}} 3^{\frac{1}{2}-2M} 2^{4M + \frac{1}{2}} P(M),
\end{equation}
where in obtaining Eq. \eqref{1/07-2} the use of Eq. \eqref{27.10.13-9} was again made.

Further transformations of Eq. \eqref{1/07-1} are necessary in order to take the full advantage of Eq. \eqref{30/06-3}. For that purpose we apply the Eq. \eqref{15/06-13} to Eq. \eqref{1/07-1}:
\begin{equation}\label{2/07-1}
G^{1, 4}_{4, 4}\left(\frac{1}{z}\Big\vert {a_{1}, \ldots, a_{4} \atop b_{1}, \ldots, b_{4}}\right) = G^{4, 1}_{4, 4}\left(z\Big\vert {1-b_{1}, \ldots, 1-b_{4} \atop 1-a_{1}, \ldots, 1-a_{4}}\right),
\end{equation}
where $(a_{p})$, $p=4$, and $(b_{q})$, $q=4$ can be read off Eq. \eqref{1/07-1}. Then $G(M, z)$ becomes
\begin{align}\label{2/07-2}
\begin{split}
G(M, z) & = r_{G}(M) \\ & \times G^{4, 1}_{4, 4}\left(-\frac{1}{Rz}\Big\vert {1, \frac{4}{3} + \frac{2M}{3}, \frac{5}{3} + \frac{2M}{3}, 2 + \frac{2M}{3} \atop \frac{3}{4} + \frac{M}{2}, 1 + \frac{M}{2}, \frac{5}{4} + \frac{M}{2}, \frac{3}{2} + \frac{M}{2}}\right),
\end{split}
\end{align}
which permits to evaluate
\begin{align}\label{2/07-3}
\begin{split}
\frac{1}{z}G&\left(M, \frac{1}{z}\right) = r_{G}(M)\, \frac{1}{z} \\ & \times G^{4, 1}_{4, 4}\left(-\frac{z}{R}\Big\vert {1, \frac{4}{3} + \frac{2M}{3}, \frac{5}{3} + \frac{2M}{3}, 2 + \frac{2M}{3} \atop \frac{3}{4} + \frac{M}{2}, 1 + \frac{M}{2}, \frac{5}{4} + \frac{M}{2}, \frac{3}{2} + \frac{M}{2}}\right). 
\end{split}
\end{align}
Eq. \eqref{2/07-3} will be transformed now using Eq. \eqref{15/06-13b}, where $\mu = -1$:
\begin{equation}\label{2/07-4}
\frac{1}{z}G\left(M, \frac{1}{z}\right) = r_{G}(M)\, G^{4, 1}_{4, 4}\left(-\frac{z}{R}\Big\vert {a'_{1}-1, \ldots, a'_{4}-1 \atop b'_{1}-1, \ldots, b'_{4}-1}\right),
\end{equation}
where now $(a'_{4})$ and $(b'_{4})$ are read off the lists in Eq. \eqref{2/07-3}; it gives finally
\begin{align}\label{2/07-5}
\begin{split}
\frac{1}{z}G&\left(M, \frac{1}{z}\right) = - \frac{r_{G}(M)}{R}\\ & \times G^{4, 1}_{4, 4}\left(-\frac{z}{R}\Big\vert {0, \frac{1}{3} + \frac{2M}{3}, \frac{2}{3} + \frac{2M}{3}, 1 + \frac{2M}{3} \atop \frac{M}{2}-\frac{1}{2}, \frac{M}{2}-\frac{1}{4}, \frac{M}{2}, \frac{M}{2}+\frac{1}{4}}\right).
\end{split}
\end{align}
It is instructive to write Eq. \eqref{30/06-3} now exclusively using the Meijer G-function in Maple notation and retaining all the multiplicative constants:
\begin{multline}\label{2/07-6}
-\frac{r_{G}(M)}{R}\,{\rm MeijerG}\big(\big[\big[0\big], \big[\ulamek{1}{3} + \ulamek{2M}{3}, \ulamek{2}{3} + \ulamek{2M}{3}, 1 + \ulamek{2M}{3}\big]\big], \\ \big[\big[\ulamek{M}{2}-\ulamek{1}{2}, \ulamek{M}{2}-\ulamek{1}{4}, \ulamek{M}{2}, \ulamek{M}{2}+\ulamek{1}{4}\big], \big[\,\,\,\big]\big], -\ulamek{z}{R}\big)\\
= r_{W}(M) \int_{0}^{R} \frac{\D x}{z - x}\\ {\rm MeijerG}\big(\big[\big[\,\,\,\big], \big[0, \ulamek{1}{3} + \ulamek{2M}{3}, \ulamek{2}{3} + \ulamek{2M}{3}, 1 + \ulamek{2M}{3}\big]\big], \qquad\qquad\\ \big[\big[\ulamek{M}{2}-\ulamek{1}{2}, \ulamek{M}{2}-\ulamek{1}{4}, \ulamek{M}{2}, \ulamek{M}{2}+\ulamek{1}{4}\big], \big[\,\,\,\big]\big], \ulamek{x}{R}\big), \qquad z > R.
\end{multline}
The same formula presented in the traditional notation, i.e.:
\begin{multline}\label{5/07-1}
\frac{r_{G}(M)}{R} G^{\,4, 1}_{4, 4}\left(-\frac{z}{R}\Big\vert {0; \frac{1}{3} + \frac{2M}{3}, \frac{2}{3} + \frac{2M}{3}, 1 + \frac{2M}{3} \atop \frac{M}{2}-\frac{1}{2}, \frac{M}{2}-\frac{1}{4}, \frac{M}{2}, \frac{M}{2}+\frac{1}{4}}\right) \\ = r_{W}(M)\,\int_{0}^{R} \frac{\D x}{x - z} \\ G^{\,4, 0}_{4, 4}\left(\frac{x}{R}\Big\vert {\,;0, \frac{1}{3} + \frac{2M}{3}, \frac{2}{3} + \frac{2M}{3}, 1 + \frac{2M}{3} \atop \frac{M}{2}-\frac{1}{2}, \frac{M}{2}-\frac{1}{4}, \frac{M}{2}, \frac{M}{2}+\frac{1}{4}}\right), \quad z > R,
\end{multline}
and in shorter notation of Eq. \eqref{27.10.13-10} Eq. \eqref{5/07-1} becomes
\begin{multline}\label{5/07-1a}
\frac{r_{G}(M)}{R} G^{\,4, 1}_{4, 4}\left(-\frac{z}{R}\Big\vert {0; \Delta(3, 2M + 1) \atop \Delta(4, 2M-2)}\right) \\= r_{W}(M) \int_{0}^{R}\!\! \frac{\D x}{x - z} G^{\,4, 0}_{4, 4}\left(\frac{x}{R}\Big\vert {\,;0, \Delta(3, 2M+1) \atop \Delta(4, 2M-2)}\right), 
\end{multline}
where $z > R$. The validity of Eq. \eqref{5/07-1} has been independently verified numerically. Eqs. \eqref{5/07-1} and \eqref{5/07-1a} appear to be less transparent than Eq. \eqref{2/07-6} and are rather more error prone. We slightly overstretched the notation of Eq. \eqref{16/06-2} in Eq. \eqref{5/07-1} by (temporarily) introducing the semicolons to explain the correct position of $0$ in coefficient lists. We stress that it is essential to keep the multiplicative constants $r_{G}(M)$ and $r_{W}(M)$ on both sides of Eqs. \eqref{2/07-6} and \eqref{5/07-1} in order to consider these equations as full solutions of Eq. \eqref{22/06-6}. The attentive reader will rapidly notice that $r_{G}(M)/r_{W}(M) = R$, and after this simplification Eqs. \eqref{2/07-6} and \eqref{5/07-1} become "bare" relations between Meijer G-functions.

For reader's convenience we quote below the formula which results from the composition of Eqs. \eqref{15/06-12} - \eqref{15/06-13b} which allows quasi-automatically to arrive at the coefficient lists appearing in $\ulamek{1}{z} G(M, \frac{1}{z})$ here, as well as serving for related problems. \\
Starting with ${_{p}F_{q}}({a_{1}, \ldots, a_{p} \atop b_{1}, \ldots, b_{q}}; Rz)$ as in Eqs. \eqref{30/06-4} and \eqref{30/06-4a}, one obtains:
\begin{multline}\label{5/07-2}
\frac{1}{z} {_{p}F_{q}}\left({a_{1}, \ldots, a_{p} \atop b_{1}, \ldots, b_{q}}; \frac{R}{z}\right) = - \frac{\prod_{k=1}^{q}\Gamma(b_{k})}{\prod_{k=1}^{p}\Gamma(a_{k})}\, \frac{1}{R}\\ \times G^{\,p, 1}_{q+1, p}\left(-\frac{z}{R}\Big\vert {0, b_{1}-1, \ldots, b_{q}-1 \atop a_{1}-1, \ldots, a_{p}-1} \right)
\end{multline}
for $p \leq q+1$, applicable in our context only for the cases when the ogf is a {\em single} generalized hypergeometric function ${_{p}F_{q}}$.

\begin{figure}[!h]
\begin{center}
\includegraphics[scale=0.22]{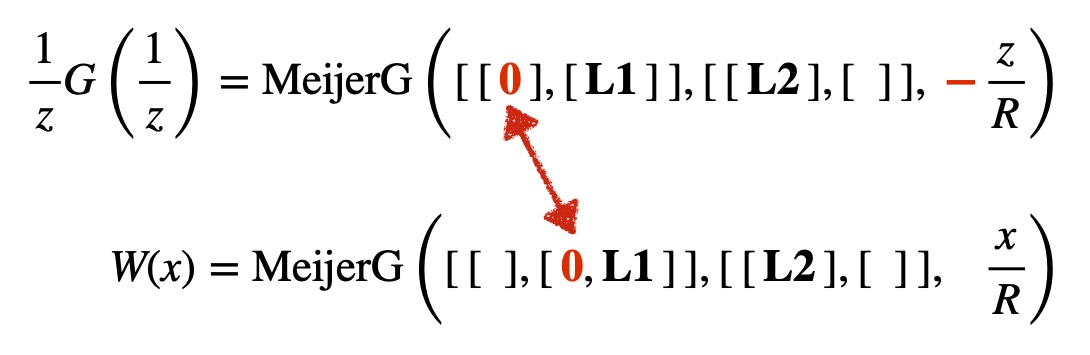}
\caption{\label{fig3}(Color online) Schematic illustration of relations of Eq. \eqref{30/06-3} for specific case of Tutte numbers of Eq. \eqref{22/06-5}, s. Eq. \eqref{2/07-6}. The lists $\boldsymbol{{\rm L1}}$ and $\boldsymbol{{\rm L2}}$ are defined before Eq. \eqref{5/07-10} in the text. We emphasize that both functions illustrated here are of Meijer G-type, but they are different functions. We neglect any multiplicative numerical constants in this illustration.
} 
\end{center}
\end{figure}

In the language of Meijer G-functions Eqs. \eqref{2/07-6} and \eqref{5/07-1} display a visibly regular scheme, which can be symbolically written down if we denote $\boldsymbol{{\rm L1}} = ( \ulamek{1}{3} + \ulamek{2M}{3}, \ulamek{2}{3} + \ulamek{2M}{3}, 1 + \ulamek{2M}{3})$ and $\boldsymbol{{\rm L2}} = (\ulamek{M}{2}-\ulamek{1}{2}, \ulamek{M}{2}-\ulamek{1}{4}, \ulamek{M}{2}, \ulamek{M}{2}+\ulamek{1}{4})$. Then, neglecting for now the multiplicative constants
\begin{equation}\label{5/07-10}
\frac{1}{z} G\left(M, \frac{1}{z}\right) \cong {\rm MeijerG}\left([\,[\,\boldsymbol{0}\,], [\,\boldsymbol{{\rm L1}}\,]\,], [\,[\,\boldsymbol{{\rm L2}}\,], [\,\,\,]\,], -\frac{z}{R}\right)
\end{equation}
and
\begin{equation}\label{5/07-11}
W_{M}(x) \cong {\rm MeijerG}\left([\,[\,\,\,], [\,\boldsymbol{0},\, \boldsymbol{{\rm L1}}\,]\,], [\,[\,\boldsymbol{{\rm L2}}\,], [\,\,\,]\,], \frac{x}{R}\right)
\end{equation}
are related through the integral formula Eq. \eqref{30/06-3} whose specific realizations are Eqs. \eqref{2/07-6} and \eqref{5/07-1}. From two previous equations we observe that reinserting the multiplicative constants one can construct the weight $W_{M}(x)$ by simply moving the number $\boldsymbol{0}$ from first bracket in Eq. \eqref{5/07-10} to second bracket in Eq. \eqref{5/07-11}, where $\boldsymbol{0}$ joins the list $\boldsymbol{{\rm L1}}$. The position of the list $\boldsymbol{{\rm L2}}$ stays unchanged in the third bracket, and the argument of $W_{M}(x)$ becomes $x/R$. Schematic display of ingredients of Eq. \eqref{30/06-3} are presented in Fig. \ref{fig3}. 

We believe that the moments $A(M, n)$ belong to a larger family of similar types of moments, for which the aforementioned reshuffling of lists gives explicitly $W(x)$ from the data of the appropriate $G(z)$, as in Eqs. \eqref{5/07-10} and \eqref{5/07-11}. If so, then there is no need to perform the inverse Mellin transform from the moments, since all the informations are already contained in the ogf $G(z)$. We are searching for possible candidates to extend the sequence $A(M, n)$ studied here. The integral relation Eq. \eqref{30/06-3} rewritten as Eq. \eqref{2/07-6} can be viewed as a variant of one-sided, finite Hilbert transform \cite{FWKing09}. However the strict condition $z > R$ imposed by convergence, requires a special care in all the manipulations.

\section{Discussion and Conclusions}\label{sec5}

We have exactly solved the moment problem of Eq. \eqref{22/06-6} following two different, and seemingly unrelated paths. The first method used was the inverse Mellin transform which resulted in exact and explicit expression for the weight functions $W_{M}(x)$ formulated in the language of Meijer G-functions. The second, less orthodox approach, consists in "upgrading" the notation for the ogf of moments $G(M, z)$, which initially was a generalized hypergeometric function ${_{4}F_{3}}$, to express it in terms of a Meijer G-function. This procedure has revealed a hitherto hidden relation between the parameter lists of $G(M, z)$ and the parameter lists of solutions $W_{M}(x)$. That observation is quite fertile, as it allows one, almost automatically, to obtain explicit forms of $W_{M}(x)$, without having recourse to any further manipulations. We believe that the sequence $A(M, n)$ of Eq. \eqref{22/06-5} belongs to a larger family of moment sequences for which analogous relations of type Eqs. \eqref{5/07-10} and \eqref{5/07-11} hold. This feature is under active consideration.

\section*{Acknowledgments}

KG thanks the LPTMC at Sorbonne Universit\'{e} for hospitality. Special thanks are due to Prof. B. Delamotte, the director of LPTMC. KG research stay at LPTMC was financed by the "Long-term research visits" program of PAN (Poland) and CNRS (France). GD wants to express his gratitude to LIPN (CNRS UMR 7030) for hosting his research. \\
KG and AH research was supported by the NCN Research Grant OPUS-12 No. UMO-2016/23/B/ST3/01714. KG acknowledges financial support by the NCN-NAWA Research Grant Preludium Bis 2 No. UMO-2020/39/O/ST2/01563.


\begin{thebibliography}{99}
\bibitem{F1}
Notice that in Mathematica notation the Meijer G-function are represented in the form 
\begin{multline*}
\text{MeijerG}[\{\{\alpha_{1}, \ldots, \alpha_{n}\}, \{\alpha_{n+1}, \ldots, \alpha_{p}\}\}, \\ \{\{\beta_{1}, \ldots, \beta_{m}\}, \{\beta_{m+1}, \ldots, \beta_{q}\}\}, x].
\end{multline*}
\end{thebibliography}
\end{document}